\documentclass[11pt,reqno]{amsart}
\usepackage{amsmath,amssymb,graphicx,amscd,amsfonts,psfrag,fancybox,qsymbols,indentfirst,latexsym}
\usepackage{calligra}
\usepackage{mathrsfs}
\usepackage{dsfont}
\usepackage[latin1]{inputenc}
 \newcommand{\ii}{{\mathrm{i}}}
\newcommand{\e}{{\mathrm{e}}}

\newcommand{\DOS}{\operatorname{DOS}}
\newcommand{\Var}{\operatorname{Var}}
\newcommand{\SC}{\operatorname{SC}}

\newcommand{\Supp}{\mathrm{Supp}\,}
\newcommand{\TW}{\mathrm{TW}_\beta\,}
\newcommand{\Lip}{\mathrm{Lip}\,}

\DeclareMathAlphabet{\mathcalligra}{T1}{calligra}{m}{n}
\newtheorem{lem}{Lemma}[section]
\newtheorem{defi}[lem]{Definition}
\newtheorem{theo}[lem]{Theorem}
\newtheorem{cor}[lem]{Corollary}
\newtheorem{prop}[lem]{Proposition}

\newtheorem{rem}[lem]{Remark}

\newtheorem{assu}[lem]{Assumption}

\theoremstyle{definition}

\theoremstyle{remark}

\def \be{\begin{eqnarray*}}
\def \ee{\end{eqnarray*}}
\def \ben{\begin{eqnarray}}
\def \een{\end{eqnarray}}

\def\sn{^{(N)}}
\def\subn{_{(N)}}

\numberwithin{equation}{section}

\begin{document}

\title[Near-extreme eigenvalues in the beta-ensembles]{Large deviations  for near-extreme eigenvalues in the beta-ensembles}
\date{\today}

\author{C. Donati-Martin}
\address{Laboratoire de Mathématiques de Versailles, UVSQ, CNRS, 
Université Paris-Saclay, 45 avenue des Etats-Unis,
F-78035 Versailles Cedex}
\email{catherine.donati-martin@uvsq.fr}
%\urladdr{http://lmv.math.cnrs.fr/annuaire/catherine-donati/}

\author{A. Rouault}
 \address{Laboratoire de Mathématiques de Versailles, UVSQ, CNRS, 
Université Paris-Saclay, 45 avenue des Etats-Unis,
F-78035 Versailles Cedex}
 \email{alain.rouault@uvsq.fr}
%\urladdr{http://rouault.perso.math.cnrs.fr/}

\subjclass[2010]{15B52, 60G57, 60F10}
\keywords{Beta ensembles, spectral distribution, top eigenvalue, large deviations}

\begin{abstract}
For beta-ensembles with convex polynomial potentials,
%For the model of log-gases on $\mathbb R$,
 we prove a large deviation principle for the empirical spectral distribution seen from the rightmost particle. 
This modified spectral distribution was introduced by  Perret and Schehr (J. Stat. Phys. 2014) to study the crowding near the maximal eigenvalue,  in the case of the GUE.  We prove also convergence of fluctuations.
\end{abstract}

\maketitle

\section{Introduction}
In random matrix models, the most popular statistics is the empirical spectral distribution (ESD). For a $N\times N$ matrix with real 
 eigenvalues $(\lambda_1, \cdots, \lambda_N)$ 
\begin{equation}\label{ESD}\mu\sn := \frac{1}{N}\sum_{k=1}^N \delta_{\lambda_{k}}\,.\end{equation}
The first step in asymptotic study is to prove the convergence of $\mu\sn$ and also of 
 the so called  integrated density of states $\mathbb E \mu\sn$. The limiting distribution $\sigma$ is more often compactly supported. A second step is  to prove the convergence of the largest eigenvalue $\lambda_{(N)} = \max (\lambda_1, \cdots, \lambda_N)$ to the end of the support of $\sigma$. At a more precise level, it is sometimes possible to establish large deviations. In the so-called $\beta$-models, the density of eigenvalues is 
\begin{equation}
\label{defpnv}\mathbb P^N_{V, \beta}(d\lambda_1, \cdots, d\lambda_N) = (Z^N_{V, \beta})^{-1}|\Delta(\lambda)|^\beta \exp \left( -\frac{N \beta}{2} \sum_1^N V(\lambda_k)\right)\ \prod_1^N d\lambda_k\end{equation} 
where $ \Delta(\lambda)$ is the Vandermonde determinant.
Under convenient assumptions on the potential $V$, the ESD satisfy the large deviation Principle (LDP) with  speed  $\beta N^2/2$ and good rate function
\begin{equation}
\label{iv}I_V (\mu) = -\Sigma(\mu) + \int V d\mu - c_V \ \ \ \hbox{if} \ \ \int |V| d\mu < \infty\,,\end{equation}
where $\Sigma$ is the logarithmic entropy
\begin{equation} \label{entropy_a} \Sigma (\mu) = \int \!\int \ln (|x-y|) d\mu(x) d\mu(y)\,, \end{equation}
and 
\begin{equation} \label{cv} c_V = \inf_\nu  -\Sigma(\nu) + \int V d\nu\,.\end{equation}
Moreover $I_V$ achieves its
 minimum $0$ 
 at a  unique probability measure $\mu_V$ which is compactly supported, and which is consequently the limit of $\mu\sn$. 

The most famous example is the Gaussian Unitary Ensemble which corresponds to $V(x) = x^2/2$ and $\beta = 2$. The limiting distribution $\mu_V$  is then the semicircle distribution~:
\begin{equation}\label{SC}\mu_{\SC}(dx) = \frac{1}{2\pi} \sqrt{4-x^2}\ 1_{[-2, 2]}(x)\, dx\,,\end{equation}
 and the result of large deviations is due to \cite{BA-G}, with $c_V= 3/4$. 

Moreover under appropriate conditions again, the support of  $\mu_V$ is an interval $[a_V,b_V]$ and  the maximal eigenvalue $\lambda_{(N)}$ converges to $b_V$. 

To analyze the "crowding" phenomenon near the largest eigenvalue, Perret and Schehr proposed in 
  \cite{perret2014near} and \cite{perret2015} to study the empirical measure~:
\begin{equation} \label{muDOS_a}
\mu_N := \frac{1}{N-1} \sum_{k=1}^{N-1} \delta_{\lambda_{(N)} - \lambda_{(k)}} \in  M_1(\mathbb R^+)\,,
\end{equation}
where $\lambda_{(1)} < \lambda_{(2)} < \cdots <\lambda_{(N)}$ are the eigenvalues of $M_N$  ranked increasingly.
They considered the Gaussian case with the Dyson values $\beta = 1,2,4$ and 
made a complete study of % the density of states near the maximum, i.e.
  $\mathbb E \mu_N$, in the limit $N \rightarrow \infty$ both in the bulk and at the edge.

In the present paper, we consider more general potentials $V$, actually convex polynomials of even degree. We first prove that $\mu_N$ converges  in probability to  the 
pushforward $\nu_V$ of $\mu_{V}$ by the mapping $x \mapsto b_V-x$. Then we prove that the family of distributions of $(\mu_N)_N$ satisfies the LDP with speed $N^2$ and a ``new" rate function which we call $I_V^{\DOS}$, referring to the name ``Density of States near the maximum" given by Perret and Schehr to $\mathbb E \mu_N$. There are two striking facts. The first one is that  the LDP is obtained for a Wasserstein topology (and not for the usual weak topology). This ensures in particular that the rate function is lower semicontinuous. The second one is that the LDP is weak i.e. we do not have a large deviation upperbound for closed sets but only for compact sets.  This implies that we could not deduce the convergence to the limit from the LDP as usual. In the Gaussian case, we have $V(x) = x^2/2$ and 
\begin{equation} \label{taux_a} 
I^{\DOS}_V(\nu) = -\Sigma(\nu) + \frac{1}{2} \Var \nu - \frac{3}{4}\,,
\end{equation}
where
for $\nu \in M_1(\mathbb R^+) $ such that $\int x d\nu(x) < \infty$, we define
\begin{equation}\label{defvar}\Var \nu = \int x^2 d\nu(x) - \left(\int x d\nu(x)\right)^2 \in [0, \infty]\,.\end{equation}

Section \ref{2} is devoted to LDPs : Proposition \ref{LDP2dim} and Corollary \ref{theoLDPV_2} study the pair $(\lambda^{(N)}, \mu_N)$, which prepares the main result, the LDP for  $(\mu_N)$ in Theorem    \ref{theoLDPV_a}. The proofs are in Section \ref{3}.  To complete the description of the asymptotic behavior of $\mu_N$ we prove also the convergence of fluctuations in Section \ref{4}. Finally, in the Appendix of Section \ref{five}, we gather some properties of the Wasserstein distance on probability measures.  
%%%%%%%%%%%%%%%%%%%%%%%%%%%%%%%%%%%%%%%%%%%%%%%%%%%%%%%%%%

\section{Assumptions and main result}
\label{2}
To begin with, let us recall the definition of the Wasserstein distance.
\begin{defi}
 Let $p \in [1, \infty[$, and $M_1^p(\mathbb R) = \{ \nu \in M_1(\mathbb R), \int |x|^p d\nu(x) < \infty \}$. For two probabilities $\mu$ and $\nu$ in $M_1^p(\mathbb R)$, the Wasserstein distance of order $p$ is defined by 
 \begin{equation}
 d_{W_p} (\mu, \nu) =\left( \inf_{ \pi \in \Pi(\mu, \nu)} \int_\mathbb R |x-y|^p d\pi(x,y) \right)^{1/p}
 \end{equation}
 where $\Pi(\mu, \nu)$ is the set of probabilities on $\mathbb R^2$ with first marginal $\mu$ and second marginal $\nu$.
 \end{defi}
Besides, we denote by $d$ the usual distance for the weak topology, given by Lipschitz bounded functions. It is known that
\begin{equation}\label{reldist}
d \leq d_{W_1} \leq d_{W_q} \ \ \hbox{for} \ q\geq 1\,.\end{equation}
 We assume that
\begin{assu}
\label{assu}
$V$ is a  convex polynomial of even degree $p \geq 2$.
\end{assu}
This assumption guarantees that $\mu_V$ is unique,  with support $[a_V, b_V]$. Moreover we have
 \begin{theo} 
 \label{theoLDP}
The sequence of distributions of $(\mu\sn)_N$
satisfies a large deviation principle  in $(M_1 (\mathbb R) , d)$,  
 with speed $\beta N^2/2$ with good rate function $I_V$ given by (\ref{iv}).
\end{theo}
This result  is Th. 2.6.1 in  \cite{agz}. As we will prove in the following section, it can be improved :
\begin{cor}
\label{cvWp}
The LDP still holds on $M_1^q$, endowed with the distance $d_{W_q}$ for any $q<p$.
\end{cor}
Fot the largest eigenvalue, we have
\newpage
\begin{prop}
\label{epsilon}
Under Assumption \ref{assu}, 
\begin {enumerate}
\item $\lambda_{(N)}$ converges in probability to $b_V$,
\item the sequence of distributions of $(\lambda_{(N)})_N$ satisfies a large deviation principle with speed $\beta N/2$, with a good rate function $J_V^+$ satisfying $J_V^+(x) = + \infty$ for $x < b_V$, i.e. 
\begin{equation}
\label{ldpright}
\lim_N \frac{2}{\beta N} \ln \mathbb P_{V, \beta}(\lambda_{(N)} > x) = - J^+_V(x)\ \ , \ \  x > b_V 
\end{equation}
with $J^+(b_V) = 0$\,,
\item the sequence of distributions of $(\lambda_{(N)})_N$  satisfies a large deviation principle with speed $\beta N^2/2$, with a rate function $J_V^-$,  on the left of $b_V$ i.e.
%, with the rate function %$J$ satisfying $J(x) = + \infty$ for $x > b_V$, i.e.
\begin{equation}
\label{ldpleft}
\lim_N \frac{2}{\beta N^2} \ln \mathbb P_{V, \beta}(\lambda_{(N)} \leq x) = - J_V^-(x) := - \inf_{\mu : \Supp \mu \subset (-\infty, x]} I_V(\mu)\ \ , \ \  x < b_V \,.
\end{equation}
\end{enumerate}
\end{prop}
Points (1) and  (2) are in \cite{agz} Prop. 2.6.6, but are more readable in \cite{Borot} Prop. 2.1. For  Point (3), see \cite{Joshape} Rem. 2.3, \cite {Feral} Sect. 4.2,  \cite{dean2006large}  and \cite{vivo2007large}).

 We are now interested in the behavior of $\mu_N$. First, we have the following convergence result: 
 \begin{prop}
\label{cvmudos}
 We denote by $\tau_c \mu$ the probability defined by
 \[\int f(x) \tau_c \mu(dx) = \int f (c-x) \mu(dx)\,.\]
 Then, as $N \rightarrow \infty$, 
$ \mu_N$ converges weakly in probability  to the probability measure~: \begin{equation}\label{nuV}\nu_V := \tau_{b_V}\mu_{V}\,.\end{equation}
 \end{prop}
Our main result rules the large deviations of the pair $(\lambda\subn, \mu_N)$. We equip  $M_1^p(\mathbb R^+)$ with $d_{W_p}$ and denote by $B(\mu ; \delta)$  the ball around $\mu$ of radius $\delta$.

We define
\begin{equation}
\label{defcalI}
\mathcal I_V (c, \nu) = I_V (\tau_c \nu)\,,
\end{equation}
which, since $\Sigma$ is invariant by the transformation $\tau_c$, is also
\begin{equation} \label{tauxV_2}
\mathcal I_V(c, \nu) = -\Sigma(\nu) + 
\int V d\tau_c \nu -c_V\,.
\end{equation}
\begin{prop} 
\label{LDP2dim}
%\label{lowerupper}
We have
\begin{enumerate}
\item
For any $c \in \mathbb R$ and $\mu \in M_1^p(\mathbb R^+)$,
\begin{eqnarray}
\nonumber
\lim_{\delta \searrow 0, \delta' \searrow 0} \liminf_{N \rightarrow \infty} \frac{2}{\beta N^2} \ln( \mathbb P_{V,\beta}^N ( \lambda\subn \in [c-\delta', c+\delta'] ,  \mu_N \in B_{W_p}(\mu, \delta)))\\ \label{lbOo} \geq - \mathcal I_V(c, \mu)\,.
\end{eqnarray}
\item
For any closed set $F\subset \mathbb R$ and $\mu \in M_1^p(\mathbb R^+)$,
\begin{eqnarray}
\nonumber
\lim_{\delta \searrow 0} \limsup_{N \rightarrow \infty} \frac{2}{\beta N^2} \ln(  \mathbb P_{V,\beta}^N  ( \lambda\subn \in F ,  \mu_N \in B_{W_p}(\mu, \delta)))\\ \label{ubF} \leq - \inf_{c \in F} \mathcal I_V(c, \mu)\,.
\end{eqnarray}
\end{enumerate}
\end{prop}

\begin{cor} \label{theoLDPV_2}
The sequence of distributions of $(\lambda\subn, \mu_N)_N$, under $\mathbb P^N_{V, \beta}$,  satisfies a weak LDP 
on $\mathbb R \times  M_1^p(\mathbb R^+)$ equipped with the product topology, 
 at speed $\beta N^2/2$ with rate function $\mathcal I_V$.
\end{cor}

From these results, we may deduce one the one hand a weak LDP for the random measure $\mu_N$, and on the other hand a conditional LDP for $\mu_N$, knowing $\lambda\subn$.
 
\begin{theo}
\label{theoLDPV_a}
 The sequence of distributions of $(\mu_N)_N$, under $\mathbb P^N_{V, \beta}$,  satisfies a weak LDP in $(M_1^p(\mathbb R^+), d_{W_p})$
 at speed $\beta N^2/2$ with rate function% $I_V^{\DOS}$.
\begin{equation} \label{tauxV_a}
I_V^{\DOS}(\nu) := \inf_{c \in (-\infty, \infty)} \mathcal I_V(c, \nu) = -\Sigma(\nu) + 
G_V(\nu) -c_V\,,
\end{equation}
with
\begin{equation}
\label{defGV}
G_V(\nu) := \inf_{c \in (-\infty, \infty)} \int V d\tau_c \nu 
\,.
\end{equation}
\end{theo}

The properties of $I_V^{\DOS}$ and $G_V$  are ruled by the following lemma:
\begin{lem}
\label{lem1}
Let $p-1 \leq q \leq p$.
\begin{enumerate}
\item
 The infimum in (\ref{defGV}) is reached at a unique point which we call $\kappa_V(\nu)$.
\item   $\nu \mapsto \kappa_V(\nu)$ is continuous for $d_{W_q}$.
\item
$\nu \mapsto G_V(\nu) = \int V(\kappa_V(\nu) -x) d\nu(x)$ is lower semicontinuous for $d_{W_q}$.
\item $I_V^{\DOS}$ is  well defined on $M_1^q(\mathbb R)$ with values in $[0, + \infty]$ and lower semicontinuous for the $W_q$ topology.
\item For  $b  \geq b_V$,  $$I_V^{\DOS}(\tau_b\mu_V) = 0.$$
\end{enumerate}
\end{lem}
%From Proposition \ref{LDP2dim}, we can deduce a LDP for $(\mu_N)_N$. 
\iffalse
\begin{theo}
\label{theoLDPV_a}
 The sequence of distributions of $(\mu_N)_N$, under $\mathbb P^N_{V, \beta}$,  satisfies a weak LDP in $(M_1^p(\mathbb R^+), d_{W_p})$
 at speed $\beta N^2/2$ with rate function $I_V^{\DOS}$.
\end{theo}
\fi
%\begin{rem}

 It follows from the property 5 in Lemma \ref{lem1}  that $I_V^{\DOS}$ is not a good rate function since the level sets are not compact.  Then, 
%Since $I_V^{\DOS}$ is not a good rate function,
  $(\mu_N)$ is not exponentially tight in scale $N^2$ (see \cite[Lemma 1.2.18]{DZ}) and we do not know if the large deviations upper bound in Theorem  \ref{theoLDPV_a}  is true for closed sets. Nevertheless, we can prove exponential tightness in a weaker topology, conditionally that $\lambda\subn$ remains bounded, which leads to :%\end{rem}

\begin{prop} Let $p-1 \leq q <p$. \\
\label{compactcond}
For  any closed set ${\mathbf F}$ of $M_1^q (\mathbb R^+)$ and any $C > b_V$, we have 
\begin{eqnarray}
\label{lbO} 
 \limsup_{N \rightarrow \infty} \frac{2}{\beta N^2} \ln( \mathbb P_{V,\beta}^N ( \mu_N \in {\mathbf F} \ | \lambda\subn \in [-C, C]))
 \leq -  \inf_{\mu \in {\mathbf F}} I_V^{\DOS}(\mu) \,.
%\mathcal I_V(c, \mu)\,.
%(- \Sigma(\mu) +  \int V d\tau_c \mu - c_V) \,.
\end{eqnarray}
\end{prop}
The conditional large deviations are ruled by the following theorem.
\newpage
%\newpage
\begin{theo} Let $p-1 \leq q < p$.
\label{theocond}
Let $\mathbf F$ be closed and $\mathbf G$ be open in $M_1^q (\mathbb R^+)$. 
\begin{enumerate}
\item If $c > b_V$, we have
\begin{eqnarray}
\label{upperdeltar}
%\limsup_{\delta'}
\limsup_N \frac{2}{\beta N^2} \ln    \mathbb P_{V,\beta}^N \left(\mu_N \in \mathbf F\ | \ \lambda\subn \in [c, c+ \delta]\right) \leq
-\inf_{\mu \in \mathbf F}\mathcal I_V^\delta(c, \mu)\\
\label{lowerdeltar}
%- \inf_{\mu \in \mathbf F, a \in [c, c+\delta]}\mathcal I_V(c, \mu)\\
\liminf_N \frac{2}{\beta N^2} \ln    \mathbb P_{V,\beta}^N \left(\mu_N \in \mathbf G\ | \ \lambda\subn \in [c, c+ \delta]\right) \geq 
%- \inf_{\mu \in \mathbf G a \in [c, c+\delta]}\mathcal I_V(c, \mu)
-\inf_{\mu \in \mathbf G}\mathcal I_V^\delta(c, \mu)\,,
\end{eqnarray}
where $\mathcal I_V^\delta(c, \mu) := \inf_{a \in [c, c+\delta]} \mathcal I_V(a, \mu)$
satisfies
\begin{equation}
\label{idelta}
\lim_{\delta \rightarrow 0}\mathcal I_V^\delta(c, \mu) = \mathcal I_V(c, \mu)\,.
\end{equation}
\item If $c < b_V$, set
\begin{eqnarray}
\label{defcalJ}
\mathcal J_V(c,\mu) &:= &\mathcal I_V(c, \mu) - \inf_{\nu\in M_1^p(\mathbb R^+)} \mathcal I_V(c, \nu)\\
 \nonumber &=& \mathcal I_V(c, \mu) - J_V^-(c),
\end{eqnarray}
(see Remark \ref{proj}).  \\
Then
\begin{eqnarray}
\label{upperdeltal}
%\limsup_{\delta'}
\limsup_N \frac{2}{\beta N^2} \ln    \mathbb P_{V,\beta}^N \left(\mu_N \in \mathbf F\ | \ \lambda\subn \in [c- \delta, c]\right) \leq - \inf_{\mu \in \mathbf F}\mathcal J_V^\delta(c, \mu)\\
\label{lowerdeltal}
%\left(\mathcal I_V(c, \mu) - \inf_{\nu\in M_1^p(\mathbb R^+)} \mathcal I_V(c, \nu)\right)\\
\liminf_N \frac{2}{\beta N^2} \ln    \mathbb P_{V,\beta}^N \left(\mu_N \in \mathbf G\ | \ \lambda\subn \in [c-\delta, c]\right) \geq - \inf_{\mu \in \mathbf G}\mathcal J^\delta_V(c, \mu)\,,
%\left(\mathcal I_V(c, \mu) - \inf_{\nu\in M_1^p(\mathbb R^+)} \mathcal I_V(c, \nu)\right)\\
\end{eqnarray}

where $\mathcal J_V^\delta(c, \mu) := %\inf_{a \in [c- \delta , c]} \mathcal I_V(a, \mu) -  \inf_{\nu\in M_1^p(\mathbb R^+)} \mathcal I_V(c, \nu)$
 \inf_{a \in [c- \delta , c]} \mathcal I_V(a, \mu) - J_V^-(c) $
satisfies
\begin{equation}
\label{jdelta}
\lim_{\delta \rightarrow 0}  \mathcal J_V^\delta(c, \mu) = \mathcal J_V(c, \mu)\,.
\end{equation}
%and $\mathcal J_V (\mu) =0$ iff $\mu = \tau_c \mu_V$.
\end{enumerate} 
\end{theo}

\medskip

\begin{rem}
\label{proj}

Let us compute the projection on $\mathbb R$ of the rate function (\ref{tauxV_2}) i.e.
\[J_V(c) := \inf_{\mu \in M_1^p(\mathbb R^+)} \mathcal I_{V}(c, \mu)\,.\]
We have, by invariance
\[J_V(c) =  \inf_{\nu : \exists \mu \in M_1^p(\mathbb R^+) : \mu = \tau_c \nu} - \Sigma(\nu) + \int V d\nu - c_V\,.\]
Recall that the support of $\mu_V$ is assumed to be $[a_V, b_V]$. 
Then, either $c \geq b_V$ and  we can take $\nu = \mu_V$, $J_V(c) = 0$, or $c < b_V$ and we get
\[J_V(c) = \inf_{\nu :\Supp \nu \subset (-\infty, c)}I_V(\nu) := J_V^-(c)\,.\]
We recover Point (3) of Prop. \ref{epsilon}. 

As  noticed in Prop. 2.1 and Rem. 2.3 in \cite{Joshape}  there is a unique $\mu$ such that $J_V(c) = I_V (\mu)$, let us call it $\mu^c$. In the Gaussian case, its explicit expression is in \cite{dean2006large} (up to some notational changes).
\end{rem}
\medskip

\begin{rem}
Let us notice that for fixed $c$, $\mathcal I_V (c, \cdot)$ and $\mathcal J_V (c, \cdot)$ may be seen as conditional rate functions. From (\ref{defcalI}), we conclude that 
\[\mathcal I_V (c, \mu) =0 \ \ \hbox{iff}\ \ \mu = \tau_c \mu_V \,,\]
whereas from (\ref{defcalJ}) and the above remark, we conclude that
\[\mathcal J_V (c, \mu) =0 \ \ \hbox{iff}\ \  \mu  = \tau_c\mu^c\,.\]
\end{rem}
\begin{rem}
\label{Ledoux}
Let us now give some additional comments relative to the Gaussian case:  $V(x) = x^2/2$. In this case, $\kappa_V(\mu)$ is the mean of $\mu$ i.e. $m(\mu) = \int x d\mu(x)$ and $G_V(\mu)$ is its variance.  
Therefore,
$$I^{\DOS}_V (\mu) = - \Sigma(\mu) + \frac{1}{2} \Var(\mu) - \frac{3}{4}\,,$$
and 
\[d\nu_V (x) = \frac{\sqrt{4x-x^2}}{2\pi} 1_{[0, 4]}(x)\!\ dx\,.\]
As we have seen above, the rate function $I^{\DOS}_V$ is zero for probabilities of the form $\mu = \tau_b \mu_{\SC}$, $b \geq 2$ and thus is not a convex rate function. Notice that this particular functional is semicontinuous not only for $d_{W_1}$ but also 
 for the weak topology. It is a consequence of  the semicontinuity of $\Var$. To prove this fact, use the representation
\[\Var(\mu) = \frac{1}{2}\!\ \mathbb E (X-Y)^2\,\]
where $X$ and $Y$ are two real random variables independent and $\mu$ distributed, and then apply Fatou's Lemma. 

A nice consequence is that the weak LDP satisfied by $\mu_N$ holds also in the weak topology (see p. 127 Remark (b) in \cite{DZ}).

\end{rem}

%%%%%%%%%%%%%%%%%%%%%%%%%%%%%%%%%%

\section{Proofs}

\label{3}
In this section, we begin with the proofs of the easiest results and we end with the proof of the main result.
\subsection{Proof of Corollary \ref{cvWp}}
The set
$$ K_M = \{ \mu \in M_1(\mathbb R), \int |V(x)| d\mu(x) \leq  M \}$$
is compact for the weak topology and is used to prove the exponential tightness for $\mu\sn$ in \cite{agz} p. 78. Actually  $K_M$ is also a compact set for the $q$-Wasserstein distance for $q<p$ (see the Appendix). It is then enough to apply  Theorem 4.2.4 of \cite{DZ}.

\subsection{Proof of Proposition \ref{cvmudos}}
 Let  $f$ a bounded Lipschitz function  with  Lipschitz constant and uniform bound less than 1. Then,
\begin{eqnarray*}
\lefteqn{ 
|\frac{1}{N-1} \sum_1^{N-1} f(\lambda_{(N)} - \lambda_{(k)}) - \frac{1}{N} \sum_1^N f(b_V- \lambda_{(k)}) |} \\
 &=& | \frac{1}{N-1}\sum_1^{N-1} (f((\lambda_{(N)} - \lambda_{(k)}) - f(b_V- \lambda_{(k)}) ) + \\
&& \qquad  \frac{1}{N(N-1)}\sum_1^{N-1} (f(b_V-\lambda_{(k)}) -  f(b_V-\lambda_{(N)})) | \\
 &\leq&  |\lambda_{(N)} -b_V| + \frac{2}{N}\,.
\end{eqnarray*}
Therefore $d( \mu_N, \tau_{b_V} \mu^{(N)}) \leq  |\lambda_{(N)} -b_V| + \frac{2}{N}$. \\
From the convergence of $\mu\sn$ to $\mu_V$, and $\lambda_{(N)}$ to $b_V$, 
 we deduce 
 that $ \mu_N$ converges to $\tau_{b_V} \mu_{V}$. $\Box$
\subsection{Proof of Lemma \ref{lem1}}
\noindent
1) Notice that the uniqueness of $\kappa_V$ comes from the convexity of $V$.\\
2)  Let $f(c) = \int V(c-x) d\nu(x)$. Then $\kappa_V(\nu)$ is the solution of \[f'(c) = \int V'(c-x) d\nu(x) = 0\,.\] Since the polynomial  $V'$ is of degree $p-1$ and, for 
$d_{W_q}$, the functions $\nu \mapsto \int x^k d\nu(x)$ are continuous for $k \leq p-1$, this implies the continuity of $\nu \mapsto \kappa_V(\nu)$.\\
3) 
Denote by $m_k(\nu)$ the $k$th moment of $\nu$. We can write $ \int V(\kappa_V(\nu) -x) d\nu(x)$ as 
$$  a_p m_p(\nu) + F(m_1(\nu), \ldots, m_{p-1}(\nu), \kappa_V(\nu))$$
where  $F$ is a polynomial function and $a_p>0$. 
The function $m_p(\nu)$ is lower semicontinuous as the supremum of the continuous functions $\int (|x|^p\wedge M) d\nu(x)$. 
The functions $\kappa_V(\nu)$ and $m_k(\nu)$, $k \leq p-1$ are continuous in $\nu$ for $d_{W_q}$. Therefore, 
$G_V$  is lower semicontinuous. \\
4)   We refer to \cite{agz} for the same properties for $I_V$, using e.g. for the positivity that $I_V^{\DOS} (\mu) = I_V( \tau_{\kappa_V(\mu)}(\mu))$. \\
From \cite{BA-G}, $-\Sigma(\mu)$ is  lower semicontinous for the topology of the weak convergence, and therefore is lower semicontinuous for the stronger topology $W_q$. \\
At last, $G_V$ is lower semicontinuous from the 3).\\
5) First notice that  $\tau_b \mu_{V}$ has a support in $\mathbb R^+$ iff $b \geq b_V$. From (\ref{tauxV_a}) and (\ref{defcalI}) we have then
\[I_V^{\DOS}(\tau_b\mu_V) = \inf_c I_V(\tau_c\tau_b \mu_V)\,,\]
and this infimum is $0$, reached at $c=b$ since $\tau_b$ is an involution. 

We could have also argued that, since the sequence $\mu_N$ converges to $\tau_{b_V} \mu_{V}$ and $I_V^{\DOS}$ is the rate function in the LDP for $\mu_N$, this insures that 
$I_V^{\DOS} (\tau_{b_V} \mu_{V}) =0$. \\

\subsection{Proof of Theorem \ref{theoLDPV_a}}

It is enough to take $c = \kappa_V(\mu)$ in the lower bound  (\ref{lbOo}) in  Proposition \ref{LDP2dim} to obtain~:
$$\lim_{\delta \searrow 0} \liminf_{N \rightarrow \infty} \frac{2}{\beta N^2} \ln( \mathbb P_{V,\beta}^N ( \mu_N \in B_{W_p}(\mu, \delta)))\\ 
 \geq -  I^{\DOS}_V( \mu)\,,$$ which implies the lower bound for open sets.\\
 For the upperbound, we take  $F = \mathbb R$ in \eqref{ubF}. $\Box$\\

\subsection{Proof of Proposition \ref{LDP2dim}}

Since the potential $V$ is assumed to be a convex polynomial, it is lower bounded by $V_{\min}$. Changing $V$ into $V - V_{\min}$ induces a change of $c_V$ into  $c_V - V_{\min}$, so in the above proofs we may and shall assume $V \geq 0$.
 
If  $\mathbb P_{V, \beta}^N$ is the distribution 
 defined in \eqref{defpnv}, we  denote by $\bar{Q}_N$ the non normalized measure
$\bar{Q}^N_{V, \beta} = Z_{V, \beta}^N \mathbb P_{V, \beta}^N$.
\iffalse
We denote by $\bar{Q}_N$ the non normalized measure
$\bar{Q}^N_{V, \beta} = Z_{V, \beta}^N \mathbb P_{V, \beta}^N$
where   $\mathbb P_{V, \beta}^N$ is the distribution 
 defined in \eqref{defpnv}.
\fi
 From \cite[p.81]{agz}
, we know that :
$$ \lim_{N \rightarrow \infty} \frac{2}{\beta N^2}\ln(Z_{V, \beta}^N) = -c_V\,. $$
Therefore, it is enough to prove  the weak LDP for the measure $\bar{Q}_N$. \\
The proof will consist in two parts : the lower bound and the upper bound.

\subsubsection{Proof of the lower bound} 
We need an approximation lemma whose  second statement is an easy consequence of  Lemma 3.3 in \cite{BA-G} (see also \cite{agz} p. 79). Indeed, the statement is given there for the distance of the weak convergence. Since the measure $\nu$ (and therefore its approximation)  has compact support, the 
same is true for the Wasserstein distance.
\begin{lem} \label{approx}

i)  Let $\mu \in M_1^p(\mathbb R^+)$, for any $\delta >0$, there exists a supported compact probability $\nu$ such that $d_{W_p}(\mu, \nu) \leq \delta$.

ii)  
 Let $\nu$ be  probability on a compact set in $\mathbb R^+$, with no atoms. 
Let $(x^{i,N})$ the sequence of real numbers defined by 
$$ x_{1,N} = \inf \{ x \geq 0 | \, \nu([0, x]) \geq \frac{1}{N}\}\,,$$
$$ x^{i+1,N} = \inf \{ x \geq x^{i,N} | \,  \nu(] x^{i+1,N}, x]) \geq \frac{1}{N}, \; 1 \leq i\leq N-2\}\,.$$
Then, 
$$x^{1,N} < x^{N-2,N} < \ldots < x^{N-1,N}\,,$$
and for any $\delta >0$ and $N$ large enough,
$$d_{W_p}\left(\nu,\frac{1}{N-1} \sum_{i=1}^{N-1} \delta_{x^{i,N}}\right) \leq \delta.$$
\end{lem}

{\bf Proof of {\it i)}} \\
 For $M>0$ and $\mu \in M_1^p(\mathbb R^+)$, we denote by $\widetilde\mu^M$ the compactly supported probability defined by $d\widetilde\mu^M = (\mu([-M,M]))^{-1} 1_{\{|x| \leq M \}} d\mu$. \\
It is easy to see that $\widetilde\mu^M$ converges weakly to $\mu$ as $M$ tends to $\infty$.
Moreover, by dominated convergence theorem,
$$ \frac{1}{ \mu([-M,M])} \int |x|^p 1_{\{|x| \leq M \}} d\mu(x) \rightarrow_{M\rightarrow \infty}  \int |x|^p d\mu(x).$$
This implies the convergence in $W_p$ distance (see Proposition \ref{propWp}, ii)). $\Box$
\medskip

To prove the lower bound (\ref{lbOo}), we will repeat almost verbatim the proof of \cite{agz} pp. 79-81, but follow step by step the rôle played by $\lambda\subn$.
We assume that $\mathcal I(c, \mu) < \infty$ so that $\mu$ has no atoms. We can also assume that $\mu$ is compactly supported, by considering $\widetilde\mu^M$ defined in Lemma \ref{approx}. One can check that $\mathcal I(c, \widetilde\mu^M) \rightarrow \mathcal I (c, \mu)$. \\
Recall that 
$$\mu_N = \frac{1}{N-1} \sum_{k=1}^{N-1} \delta_{\lambda_{(N)} - \lambda_{(k)}}$$
where the $\lambda_{(k)}$ are the increasing sequence of eigenvalues.

Then,  if $C = [c-\delta', c+\delta']$ and $B := B_{W_p}(\mu, \delta)$,
\begin{eqnarray}
\lefteqn{ 
\bar{Q}_N (\lambda\subn \in C, \mu_N  \in B)
} \nonumber\\
& = &N! \int_{ \Delta_N \cap \{\frac{1}{N-1} \sum_{k=1}^{N-1} \delta_{\lambda_N - \lambda_k} \in B , \lambda_N \in C\} }\widehat L_N(\lambda_1, \ldots \lambda_N) d\lambda_1 \ldots d\lambda_N  \label{eq1}
\end{eqnarray}
where
\begin{equation}
\label{DeltaN}
\Delta_N = \{ \lambda_1 < \lambda_2 < \ldots < \lambda_{N}\}\,,
\end{equation}
and 
\begin{eqnarray}
\lefteqn{\widehat L_N( \lambda) = \prod_{i<j} |\lambda_i -  \lambda_j|^\beta \exp( - \frac{\beta N}{2} \sum_{i=1}^N  V(\lambda_i))} \label{eq2}\\
&&=  \prod_{i<j<N-1} |\lambda'_i -  \lambda'_j|^\beta \prod_{i<N} |\lambda'_i|^\beta\exp( - \frac{\beta N}{2}( \sum_{i=1}^{N-1} V( \lambda_N -  \lambda'_i ) +  V(\lambda_N)) ) \nonumber\\
&&:= L_N((\lambda'_i)_{i<N}, \lambda_N) \nonumber
\end{eqnarray}

where $\{\lambda'_i = \lambda_N - \lambda_i, \, i<N\} $ is a decreasing family of positive numbers.
Since the density $L_N((\lambda'_i)_{i<N}, \lambda_N)$ is symmetric in $(\lambda'_i)$, we can write :
\begin{eqnarray*}
\lefteqn{ 
\bar{Q}_N (\lambda\subn \in C, \mu_N \in B) } \\
& = &N \int_{ \mathbb R_+^{N-1} \times C \, \cap \{\frac{1}{N-1} \sum_{k=1}^{N-1} \delta_{ \lambda'_k} \in B\} } L_N(\lambda'_1, \ldots , \lambda'_{N-1}, \lambda_N)) d\lambda'_1 \ldots d\lambda'_{N-1} d \lambda_N\,.
\end{eqnarray*}
From  Lemma \ref{approx}, for $N \geq N_\delta$, 
$$   \left\{ (\lambda'_k)_k : \  |\lambda'_k -x^{k,N} | \leq \frac{\delta}{2} , \, \forall k < N \right\} 
\subset  \left\{ (\lambda'_k)_k  : \  \frac{1}{N-1} \sum_{k=1}^{N-1} \delta_{ \lambda'_k}  \in B\right\} \,,$$
and we can write $L_N$ as 
\begin{eqnarray*}
L_N &=&  \prod_{i<j<N} |(\lambda'_i -x^{i,N}) -  (\lambda'_j- x^{j,N}) + x^{i,N} - x^{j,N}|^\beta \prod_{i < N} |\lambda'_i|^\beta  \\
&\times& \; \exp( - \frac{\beta N}{2}( \sum_{i=1}^{N-1}V (  \lambda_N-  x^{i,N}  - (\lambda'_i- x^{i,N}) ) +  V(\lambda_N - c +c)) )\,.
\end{eqnarray*}

Set $y_i = \lambda'_i - x^{i,N}$, $i <N$ and $y_N= \lambda_N-c$.
Then,
\begin{eqnarray*}
\lefteqn{ 
\bar{Q} (\mu_N \in B) } \\
&\geq & N \int_{\Delta_N (\delta)}
 \prod_{i<j<N} |y_i -  y_j+ x^{i,N} - x^{j,N}|^\beta \prod_{i \leq  N-1} |y_i+ x^{i,N}|^\beta  \\
&& \quad \exp( - \frac{\beta N}{2}( \sum_{i=1}^{N-1} (V ((y_N+c) - (y_i + x^{i,N}) ) + V (y_N +c)) ) \prod_{i\leq N} dy_i
\end{eqnarray*}
where 
\begin{eqnarray*}
\Delta'_N &=& \{ y_1 <  y_2 < \ldots < y_{N-1}\}\\
\Delta_N (\delta) &=& \{(y_1, \dots, y_N) : (y_i)_{i<N} \in [0, 
\delta/2]^{N-1} ,   y_N \in [-\delta, \delta]\} \cap \Delta'_N\,.
\end{eqnarray*}

Since on $\Delta'_N$, the $(y_i)$ and the $(x^{i,N})$ form both increasing sequences, we have the lower bound:
$$ |y_i -  y_j+ x^{i,N} - x^{j,N}| \geq \sup\{ | x^{i,N} - x^{j,N}|, |y_i -  y_j|\}$$
and we use the same minoration as in \cite{agz} for the term 
\begin{eqnarray*}
A&:=& \prod_{i<j<N} |y_i -  y_j+ x^{i,N} - x^{j,N}|^\beta \\
&\geq&  \prod_{i+1<j<N} | x^{i,N} - x^{j,N}|^\beta \prod_{i<N-1} | x^{i,N} - x^{i+1,N}|^{\beta/2} \prod_{i<N-1} | y_i - y_{i+1}|^{\beta/2}\,.
\end{eqnarray*}

For the second term, we use, since the $y_i$ and $x^{i,N}$ are positive,
$$ \prod_{i < N-1} |y_i+ x^{i,N}|^\beta  \geq \prod_{i < N-1} |y_i|^\beta\,.$$
We get:
\begin{equation*}
\bar{Q} (\lambda\subn \in C, \mu_N \in B)  \geq P_{N,1}P_{N, 2}\,,
\end{equation*}
where
\begin{eqnarray}
\nonumber
P_{N,1} 
&=& N\exp -\frac{\beta N}{2}\left(\sum_{i=1}^{N-1} V(c -x^{i,N}) + V(c)\right)\\
\label{prod}
&&\times \prod_{i+1 < j<N} |x^{i,N}- x^{j,N}|^\beta \prod |x^{i,N} - x^{i+1, N}|^{\beta/2}
\end{eqnarray}
and
\begin{eqnarray*}
P_{N,2} 
&=&  \int_{\Delta_N(\delta)}
\prod_{i < N-1}|y_i - y_{i+1}|^{\beta/2}|y_i|^\beta\\
&&\ \ \times \exp -\frac{\beta N}{2} \sum_{i=1}^{N-1} V(c-  x^{i,N} -y_i + y_N) - V(c- x^{i,N})\\
&&\ \ \times \exp -\frac{\beta N}{2}  \left[ V(y_N +c) - V(c)\right] \ \prod_{i\leq N} dy_i\,.
\end{eqnarray*}

Since we have assumed that $\mu$ is compactly supported, the sets $\{x^{i,N}, 1 \leq i \leq N-1\}$ are uniformly bounded and by continuity of $V$,
\[\lim_{\delta \rightarrow 0}\sup_N \sup _{1 \leq i \leq N}\sup_{|x|\leq \delta} |V(c-x^{i,N} +x)- V(c-x^{i,N}))| = 0\,,\]
and
\[\lim_{\delta \rightarrow 0}\sup_{|x|\leq \delta} |V(c+x) - V(c)| =0\,.\]
Moreover, writing $u_1 = y_1, u_{i+1} = y_{i+1}- y_i$, with $\delta'' = \min(\delta/2, \delta')$
\begin{eqnarray*}
\lefteqn{ \int_{\{(y_i)_{i<N} \in [0, \frac{\delta}{2}]^{N-1}\}  \cap \Delta'_N \cap \{ y_N \in [-\delta, \delta] \}} 
 \prod_{i < N-1}|y_i - y_{i+1}|^{\beta/2}|y_i|^\beta\ \prod_{i\leq N} dy_i}\\
&\geq&\int_{\{(u_i)_{i < N} \in [0, \frac{\delta''}{N}]^{N}} u_1^{\beta} (u_2\cdots u_{N-1})^{3\beta/2} u_N^{\beta/2}\prod_{i \leq N} du_i\\
&\geq& C_1^{-N} \left(\frac{\delta''}{N}\right)^{C_1N}
\end{eqnarray*}
for some constant $C_1$, which yields
\[\lim_{\delta, \delta' \rightarrow 0}\liminf_N \frac{2}{\beta N^2} \log P_{N,2} \geq 0\,.\]

On the other hand, from the choice of the $x^{i,N}$, we have 
$$\lim \frac{1}{N-1} \sum_{i=1}^{N-1} V(c- x^{i,N}) = \int V(c-x) d\mu(x)\,.$$
Finally the product in (\ref{prod}) can be managed exactly as in \cite{agz} p. 80.
We conclude
\begin{eqnarray*}
\lefteqn{\lim_{\delta' \searrow 0} \lim_{\delta \searrow 0} \liminf_N \frac{2}{\beta N^2} \ln( \bar{Q} (\lambda\subn \in C, \mu_N \in B) } \\
&\geq &
- \int\!\int \ln(y-x) d\mu(x) d\mu(y) -  \int  V(c-x) d\mu(x) 
\end{eqnarray*}
which is the expected lower bound. $\Box$

\subsubsection{Proof of the upper bound}
We start as in the proof of the lower bound with the representation \eqref{eq1}. Formula \eqref{eq2} can be rewritten as 
\begin{eqnarray}
\frac{2}{\beta} \ln(\widehat L_N (\lambda))  &=& 2 (N-1)^2 \int\int_{x<y} \ln(|x-y|) d\mu_N(x) d\mu_N(y) \nonumber\\
 &+& 2(N-1) \int\ln(|x|) d\mu_N(x)    \nonumber\\
 &-&   (N-1)^2 \int V(\lambda_N-x) d\mu_N(x) - (N-1) V(\lambda_N)    \nonumber \\
 &-&  \sum_{i=1}^N V(\lambda_i) \label{densite}
 \end{eqnarray}
 where $\mu_N =  \frac{1}{N-1} \sum_{k=1}^{N-1} \delta_{\lambda_N - \lambda_k} $, since we are on the set $\Delta_N$ defined in (\ref{DeltaN}). \\
 Under $\bar{Q}_N$, the $\lambda_i$ are a.s. distinct so $(\mu_N )^{\otimes 2} (\{(x,y) ; x=y\}) =\frac{1}{N-1}$ a.s.. Therefore, for every  $M \in \mathbb R$ we can write,
 \begin{eqnarray}
 \lefteqn{ - 2 \int\int_{x<y} \ln(|x-y|) d\mu_N(x)  d\mu_N(y) =} \nonumber\\
 && - 2 \int \int_{x<y} \ln(|x-y|) d\mu_N(x) d\mu_N(y) + M\left ((\mu_N )^{\otimes 2} (\{(x,y) ; x=y\}) - \frac{1}{N-1}\right) \nonumber\\
 & \geq& \int \int_{\mathbb R^2} (-\ln(|x-y|) \wedge M) d\mu_N(x) d\mu_N(y) - \frac{M}{N-1} \nonumber\\& 
:=&  - \Sigma^M(\mu) - \frac{M}{N-1}\,. \label{lnL}
\end{eqnarray}
We have assumed $V \geq 0$ so  the  second term in the third line of \eqref{densite} is non positive. 
 For the first term of the same line, notice that 
\begin{eqnarray*}
-  (N-1)^2 \int V (\lambda_N-x) d\mu_N(x) 
  &\leq&  - (N-1)^2 \inf_{c \in F} \int V(c-x) d\mu_N(x)\,.\\
 \end{eqnarray*}

 We bound the  term in the second line of \eqref{densite} by 
 $ (N-1) \int |x| d\mu_N(x).$
On the event $\{\mu_N \in B\}$, we have
 \begin{eqnarray*}
\frac{2}{\beta} \log \widehat L_N(\lambda) &\leq & - (N-1)^2 \inf_{\nu \in B} \left(- \Sigma^M(\nu) + \inf_{c \in F} \int V d\tau_c \nu\right)\\
&& + (N-1) \sup_{\nu \in B} \int |x| d\nu(x) + (N-1)M   \\
 && \quad - \sum_{i=1}^N V(\lambda_i)\,. 
 \end{eqnarray*} 

 Since
\[N! \int_{\Delta_{N}} \exp -\frac{\beta}{2} \sum_{i=1}^{N} V(\lambda_i) = \left[\int \exp -\frac{\beta}{2} V(\lambda) d\lambda\right]^{N}\,,\]
and $\sup_{\nu \in B} \int |x| d\nu(x) < \infty$,
 we obtain :
\begin{eqnarray*}\limsup_N \frac{2}{\beta N^2} \ln(  \bar{Q} (\lambda\subn \in F, \mu_N \in B) \leq
% &\leq  -  \inf_{c \in F , \nu \in B} (-\Sigma^M(\nu) + \int V d\tau_c \nu )\\
%&=
- \inf_{\nu \in B} \left(-\Sigma^M(\nu) + \inf_{c \in F}\int Vd\tau_c\nu\right)\,.
\end{eqnarray*}
To let $\delta \rightarrow 0$, we need semicontinuity in $\nu$. We know that $\Sigma^M$  is  lower semicontinuous. 
Assume that $F = [a, \infty)$. Since $V$ is convex, the infimum 
  $\inf_{c \geq a} \int V(c-x) d\nu(x) $ is reached at $c= \kappa_V(\nu)$ or at $c=a$, so that
\[ \inf_{c \geq a} \int V(c-x) d\nu(x) = \int V(\max(a, \kappa_V(\nu))-x) d\nu(x)\]
which is  a lower semicontinuous function of $\nu$. \\
We obtain : 
\begin{eqnarray}
\nonumber
 \lefteqn{\lim_{\delta \searrow 0} \limsup_N \frac{2}{\beta N^2} \ln(  \bar{Q} (\lambda^{(N)} \in F , \mu_N \in B) )\leq}\\   \label{ubF1}
&&-   \left(-\Sigma^M(\mu) +    \inf_{c \in F} \int V(c-x) d\mu(x)\right)\,.
\end{eqnarray}
and since $\Sigma^M$ grows to $\Sigma$ as $M$ goes to infinity, this yields the upper bound \eqref{ubF}.  \\
The same is true for $F = ]-\infty, a]$. \\
Now,  take $F$ a non empty closed set. If $\kappa_V(\mu) \in F$,  \eqref{ubF1} is clearly true. \\
If $\kappa_V(\mu) \notin F$, then $\kappa_V(\nu) \notin F$ for $\nu \in B$ for small $\delta$. \\
Denote by $a_- = \sup\{ x < \kappa_V(\mu), x \in F\}$ and $a_+ = \inf\{ x> \kappa_V(\mu), x  \in F\}$. Then,
$ F \subset ]-\infty, a_-] \cup [a_+, \infty[$ and 
\begin{eqnarray}\nonumber
\lefteqn{ \lim_{\delta \searrow 0} \limsup_N \frac{1}{N^2} \ln(  \bar{Q} (\mu_N \in B, \lambda\subn \in F) ) \leq} \\
&&-   \left(-\Sigma^M(\mu) +   \int V(a_- -x) d\mu(x) \wedge  \int V(a_+ -x) d\mu(x) \right)\,.\end{eqnarray}
The last term in the above equation is $ \inf_{c \in F} \int V(c-x) d\mu(x)$. $\Box$

\subsection{Proof of Propostion \ref{compactcond} and Theorem \ref{theocond}}
%Lemma \ref{quasiexptight} and Proposition \ref{addendum}} 
From Corollary \ref{theoLDPV_2}, we have
\begin{equation}\label{bi}
\limsup \frac{2}{\beta N^2} \ln \mathbb P_{V, \beta}^N( \mu_N \in \mathbf F, \lambda\subn \in [a,b]) \leq - \inf_{\mu \in \mathbf F, c \in [a,b]} \mathcal I_V(c, \mu)\,,\end{equation}
as soon as $\mathbf F$ is compact. To extend this property to closed sets, we follow the classical way and prove 
\begin{lem}
\label{quasiexptight}
Let $q<p$. For any $-\infty < a < b < \infty$ and $M > 0$, there exists a compact set $\mathbf K_{a,b, M}$ of $M_1^q(\mathbb R^+)$ such that 
\begin{equation}
\label{bigK} \limsup_N \frac{2}{\beta N^2} \ln( \mathbb P_{V,\beta}^N ( \lambda\subn \in [a,b] ,  \mu_N \notin {\mathbf K_{a,b,M}})) \leq -M\,.\end{equation}
\end{lem}
\proof
Let 
\[\mathbf K_M :=\{\mu \in M_1^p (\mathbb R^+) : \int |V| d\mu \leq M\}\,,\]
which is compact of $M_1^q(\mathbb R^+)$ from Proposition \ref{propWp}.

With our assumptions on the potential $V$, there exists $c_1,c_2 >0$  such that
\[|V(x)| \leq c_1 x^p +c_2\,,\]
Let $a<b$ and $C = \sup\{|a|, |b|\}$. \\
For $N \geq 2$, using the convexity of $x^p$
\[\{\mu_N(V) \geq M\} \cap \{\lambda_{(N)} \in [a, b]\} \subset \{\mu\sn(V) \geq M'\}\]
where
\[M = c_1(C^p 2^{p-1} + 2^{p-1} M') + c_2 \,.\]

It remains to use  the  exponential tightness for the ESD  $\mu^{(N)}$, see \cite{agz}, p. 77) where it is shown that :
$$ \limsup_N \frac{2}{\beta N^2} \ln( \mathbb P_{V,\beta}^N ( \mu^{(N)} \notin \mathbf K_{\bar{M}})) \leq -M $$
where $\bar{M}$ is an affine function of $M$. 
From Lemma \ref{quasiexptight}, \eqref{bi} is satisfied for $\bf F$ a closed set of $M_1^q(\mathbb R^+)$. $\Box$

\subsubsection{Proof of  Proposition \ref{compactcond}}
By Proposition  \ref{epsilon}, we know that for $\Delta= [-C,  C]$ large enough, then $\mathbb P_{V, \beta}^N ( \lambda\subn \in \Delta) \rightarrow 1$, so that
\[\lim_N \frac{2}{\beta N^2} \ln \mathbb P_{V, \beta}^N(\lambda\subn \in \Delta) = 0\,,\] and then, for a closed set $\mathbf F$, 
\begin{eqnarray}\nonumber
\limsup_N \frac{2}{\beta N^2} \ln \mathbb P_{V, \beta}^N( \mu_N \in \mathbf F \ | \ \lambda\subn \in \Delta) =\\
\nonumber
\limsup_N \frac{2}{\beta N^2} \ln \mathbb P_{V, \beta}^N( \mu_N \in \mathbf F , \lambda\subn \in \Delta)\\
\leq - \inf_{\mu \in \mathbf F, c \in \Delta} \mathcal I_V(c, \mu)\,,
\end{eqnarray}
from (\ref{bi}) for closed sets. Now, we use the easy bound
\[ \inf_{\mu \in \mathbf F, c \in \Delta} \mathcal I_V(c, \mu) \geq  \inf_{\mu \in \mathbf F, c \in (-\infty , \infty)} \mathcal I_V(c, \mu)= 
\inf_{\mu \in \mathbf F}  I^{\DOS}_V(\mu)\,.\]
\subsubsection{Proof of Theorem \ref{theocond}}
We  use Proposition \ref{epsilon} %and (\ref{ldpleft}) 
 to  estimate the probabilities of the  conditioning events. On the one hand (\ref{bi}) for closed sets  and   (\ref{ldpright}) lead to (\ref{upperdeltar}) and on the other hand  (\ref{bi}), (\ref{ldpleft}) and Remark \ref{proj} lead to (\ref{upperdeltal}), using that $J_V^-$ is decreasing on $[-\infty, b_V]$.

For the lower bounds, we use the lower bound coming from  the LDP for both variables (Theorem \ref{LDP2dim} (1) or Corollary \ref{theoLDPV_2}) and 
(\ref{ldpright}) and (\ref{ldpleft}), respectively.

It remains to prove (\ref{idelta}) and (\ref{jdelta}), but it is straightforward since
\begin{eqnarray*}
\lim_{\delta \rightarrow 0} \inf_{a \in [c, c+\delta]} \int V(a-x) d\mu(x) = \lim_{\delta \rightarrow 0}  \inf_{a \in [c-\delta, c]} \int V(a-x) d\mu(x)\\ = \int V(c-x) d\mu(x)\,.
\end{eqnarray*}

\section{Fluctuations}
\label{4}
We want to study the fluctuations of $\mu_N$ around  its limit $\nu_V$ given in (\ref{nuV}).  
There are two contributions: the fluctuations  of the largest eigenvalue and  the fluctuations of the ESD. This yields  a dichotomy according to the behavior of the test function.  For the sake of simplicity, we choose a simple assumption on the test function $f$ which is far from optimal. For $V$ and $\beta$ we introduce a new assumption :
\begin{assu}
\label{assunew}
$V$ satisfies Assumption \ref{assu} and $\beta =1,2,4$, or $V(x) = x^2/2$ and $\beta > 0$.
\end{assu}
\begin{prop}
Let $f$ be a bounded $\mathcal C^2$ function with two bounded derivatives. 
\begin{enumerate}
\item If $V$ and $\beta$ satisfy Assumption \ref{assunew} and if 
 $\nu_V(f') \not= 0$, 
\[N^{2/3} \left( \mu_N (f) - \nu_V (f)\right) \Rightarrow \nu_V (f') \TW\,.\]
where $TW_\beta$ denotes the Tracy-Widom distribution of index $\beta$ (see \cite{ramirez2011beta} for a definition), and where $\Rightarrow$ denotes the convergence in distribution.
\item  If $V$ satisfies Assumption \ref{assu} and  $\beta >0$  and if 
 $\nu_V(f') = 0$, 
$$ N \left( \mu_N (f) - \nu_V (f)\right) \Rightarrow \mathcal N(-f(0) + m_V(f), \sigma_V^2(f))$$
where 
\begin{equation}
\label{varfluc}\sigma_V^2(f) = \frac{1}{4\beta} \sum_{k=1}^\infty ka_k^2
\,,\end{equation}
with
\[a_k = \frac{2}{\pi} \int_0^\pi f\left(\frac{b_V - a_V}{2}(1-\cos \theta)\right) \cos k\theta\!\ d\theta\,,\]
and 
\begin{equation}
\label{bias}m_V(f)= \left(\frac{2}{\beta} -1\right)\int f(b_V -t)  d\gamma_V(t)
\end{equation}
where $\gamma_V$ is a signed  measure on $[a_V, b_V]$ given by formula (3.54) in \cite{johansson1998fluctuations}. 
\end{enumerate}
\end{prop}
Let us notice, from Remark 3.5 in \cite{johansson1998fluctuations}, that in the Gaussian case, $V(x) = x^2/2$, then
\[d\gamma_V (t) = \frac{1}{4}\delta_{-2} +  \frac{1}{4}\delta_{2}- \frac{1}{2\pi} \frac{dt}{\sqrt{4-t^2}}\,.\]

\proof
Let  $H > b_V- a_V$. Set $K_H$ the random set defined by 
%\[K_H = \{ (\lambda_{(1)}, \cdots , \lambda_{(N)}) :  \forall i,  \  -H \leq \lambda_{(N)} - \lambda_{(i)} \leq H \mbox{ and }  -H \leq b_V - \lambda_{(i)} \leq H \}\,,\]
\[K_H = \{  (\lambda_{(1)}, \cdots , \lambda_{(N)}) : b_V - H \leq \lambda_{(1)}\ ;\ \lambda_{(N)} \leq b_V + H\}\,,\]
and
\[M =   \max\left\{ \frac{1}{2} \sup_{x\in [-2H , 2H]} |f''(x)|, \sup_{x\in [-H , H]} |f(x)|,   \sup_{x\in [-H , H]} |xf'(x)|\right\}\,.\]

Setting
\begin{eqnarray*} S_N (f) := (N-1) \mu_N (f) 
=\sum_1^{N-1} f(\lambda_{(N)} - \lambda_{(k)})\,, 
\end{eqnarray*}
we make  a Taylor expansion of $f$~:
\begin{eqnarray*}
f(\lambda_{(N)} - \lambda_{(k)}) &=&  f(b_V - \lambda_{(k)} ) +   \varepsilon_N  f'(b_V - \lambda_{(k)} )+ r_{k, N}(f)\,,\end{eqnarray*}
with
 \[\varepsilon_N := \lambda_{(N)} -b_V\]
and
\[|r_{k,N}(f)|1_{K_H} \leq  M \varepsilon^2_N\,.\]
Adding,
\begin{eqnarray}
\nonumber
 S_N (f)&=& \sum_1^{N-1} f(b_V - \lambda_{(k)} ) + \varepsilon_N \sum_1^{N-1} f'(b_V - \lambda_{(k)} )+ \sum_1^{N-1}r_{k, N}\\
\label{snf}
&=& \sum_{i=1}^{N} f(b_V - \lambda_i ) + \varepsilon_N \sum_{i=1}^{N} f'(b_V - \lambda_i ) + R_N(f)
\end{eqnarray}
where 
\[R_N (f):= \sum_1^{N-1}r_{k, N}(f)\ - f(-\varepsilon_N) - \varepsilon_N f'(-\varepsilon_N)\,,\]
 satisfies
\begin{equation}
\label{RN}
|R_N(f)|1_{K_H} \leq M(N \varepsilon_N^2 + |\varepsilon_N| + 1)\,.
\end{equation} 
Setting
\begin{eqnarray}\nonumber\Delta_N(f) &:=& \sum_{i=1}^{N} f(b_V - \lambda_i )  - N\nu_V(f)\\ \label{Delta}
&=& \sum_{i=1}^{N} f(b_V - \lambda_i )  - N\int f(b_V - x) d\mu_V (x)
\,,\end{eqnarray}
(\ref{snf}) gives
\begin{eqnarray}
\label{27}
S_N(f) - N \nu_V(f) =   N\varepsilon_N\nu_V(f') +  \Delta_N (f) + \varepsilon_N \Delta_N (f') + R_N(f)\,.\end{eqnarray}
The two sources of fluctuations are the convergences of $\varepsilon_N$ (rescaled) and $\Delta_N(f)$.\\
On the one hand we  know (Prop.  \ref{epsilon}) that 
\begin{equation}
\label{largest}\varepsilon_N \rightarrow 0\ \ \hbox{in probability}\,,
\end{equation} and   the fluctuations are ruled by
\begin{equation}
\label{TW}
N^{2/3}\varepsilon_N \Rightarrow \TW\,.\end{equation}
(see \cite{deift2007universality} in the cases $\beta=1,2,4$, and \cite{ramirez2011beta} for the Gaussian case and $\beta > 0$).
\\ On the other hand, under our assumptions on $V$ and $f$,  
\begin{equation}
\label{limdelta} \Delta_N(f) 
 \Rightarrow \mathcal N (m_V(f) ; \sigma_V^2(f))\,,\end{equation}
where $m_V(f)$ and $\sigma^2(f)$ are  given by (\ref{bias}) and (\ref{varfluc}), respectively (see (\cite{johansson1998fluctuations} Theorem 2.4)).
\smallskip
\\
 1. If $\nu_V(f') \not= 0$,  we  set,  for the sake of simplicity
\[ N^{-1/3}[S_N (f) - N \nu(f)] = N^{2/3}\varepsilon_N \nu_V(f') +  R'_N (f)\,.\]
We have then, if $\Phi(f) := \mathbb E [\exp[ \ii \nu(f') \TW]$, 
\begin{eqnarray}
\nonumber
\lefteqn{\mathbb E \e^{\ii N^{-1/3}[S_N (f) - N\nu(f)]} - \Phi(f) =}\\
\nonumber
&& \mathbb E \e^{\ii N^{2/3}\varepsilon_N \nu(f')} - \Phi(f)\\ 
\nonumber
&&\ + \ \mathbb E \left(\e^{\ii N^{2/3}\varepsilon_N \nu (f')}\left[\e^{\ii R'_N(f)} - 1\right]1_{K_H}\right)\\
&&\ \ +  \ \mathbb E \left(\e^{\ii N^{2/3}\varepsilon_N \nu (f')}\left[\e^{\ii R'_N(f)} - 1\right]1_{K_H^c}\right)\,.
\end{eqnarray}
\begin{itemize}
\item
The first term converges to zero, thanks to (\ref{TW}). 
\item The second term is bounded by 
$\mathbb E\left( |R'_N(f) \wedge 2| 1_{K_H}\right)$ 
which tends to zero since on $K_H$ 
\begin{eqnarray*}|R'_N (f)|&\leq& N^{-1/3}|\Delta_N (f)| + N^{-1/3}\varepsilon_N \Delta_N (f') \\
&&+ M N^{2/3}(\varepsilon_N)^2 + MN^{-1/3}(|\varepsilon_N| + 1)\end{eqnarray*}
and each of these terms tends to zero in probability, thanks  to (\ref{limdelta}) and (\ref{largest}).
\item The third one is bounded by $2 \mathbb P((K_H)^c)$ which tends to zero, since the extreme eigenvalues tend to the endpoints of the support.
\end{itemize}
 This allows to conclude that
\[N^{2/3} \left( \mu_N (f) - \tau_{b_V}\mu_V (f)\right)  \Rightarrow \nu_V (f') \TW\,.\]
\smallskip
\\ 2. If $\nu(f') = 0$, then,
 \begin{equation}
S_N(f) - N \nu(f) =     \Delta_N (f) - f(0) + \varepsilon_N \Delta_N (f') + \tilde{R}_N(f)\,,\end{equation}
with
$$  \tilde{R}_N(f)  := \sum_1^{N-1}r_{k, N}(f)\ - (f(-\varepsilon_N) - f(0))  - \varepsilon_N f'(-\varepsilon_N)\,.$$
From (\ref{limdelta}), \[ \Delta_N (f) - f(0) \Rightarrow \mathcal N(-f(0) + m_V(f) ; \sigma_V^2(f))\,.\] 
Moreover
$\varepsilon_N \Delta_N (f') $ and $\tilde{R}_N(f) 1_{K_H}$ tend to 0 in probability.
% and $\tilde{R}_N(f) 1_{K_N} \rightarrow 0$. 
The rest of the proof goes as before.
 $\Box$

\section{Appendix}
\label{five}
 We give some properties of the Wasserstein distance $d_{W_p}$.
 \begin{defi}
 Let $p \in [1, \infty[$, and $M_1^p(\mathbb R) = \{ \nu \in M_1(\mathbb R), \int |x|^p d\nu(x) < \infty \}$. For two probabilities $\mu$ and $\nu$ in $M_1^p(\mathbb R)$, the Wasserstein distance of order $p$ is defined by 
 \begin{equation}
 d_{W_p} (\mu, \nu) =\left( \inf_{ \pi \in \Pi(\mu, \nu)} \int_\mathbb R |x-y|^p d\pi(x,y) \right)^{1/p}
 \end{equation}
 where $\Pi(\mu, \nu)$ is the set of probabilities on $\mathbb R^2$ with first marginal $\mu$ and second marginal $\nu$.
 \end{defi}
 \begin{rem} For the Wasserstein distance of order 1, we have the duality formula 
 $$ d_{W_1} (\mu, \nu) = \sup_{\Vert f\Vert_{\Lip} \leq 1} \left( \int f d\mu - \int f d\nu \right).$$
 \end{rem}
 We now give a  characterization of the convergence of probabilities in the topology induced by $d_{W_p}$ on $M_1^p(\mathbb R)$. We refer to \cite[Def. 6.8 and Theorem 6.9]{V}.\\
  In the following, we denote by $\mu_n \rightarrow \mu$ the weak convergence of probabilities, i.e. against bounded continuous functions.
 \begin{prop} \label{propWp}
 Let $(\mu_n)_{n \geq 0}$ a sequence of probabilities in $M_1^p(\mathbb R)$ and $\mu \in M_1^p(\mathbb R)$. The following assertions are equivalent :
 \begin{itemize}
 \item[\it i)] $\displaystyle d_{W_p} (\mu_n, \mu) \rightarrow 0,$
 \item[\it ii)] $\mu_n \rightarrow \mu$ and $\displaystyle  \int |x|^p d\mu_n(x) \rightarrow  \int |x|^p d\mu(x)$, 
 \item[\it iii)] $\mu_n \rightarrow \mu$ and $\displaystyle \limsup_{n\rightarrow \infty} \int |x|^p d\mu_n(x) \leq \int |x|^p d\mu(x)$, 
  \item[\it iv)] $\mu_n \rightarrow \mu$ and $\displaystyle \lim_{R \rightarrow \infty} \limsup_{n\rightarrow \infty} \int_{ |x| \geq R} |x|^p d\mu_n(x) =0$, 
  \item[\it v)] For all continuous functions $f$ with $|f(x)| \leq C(1 + |x|^p)$, one has 
  $$ \int f(x) d\mu_n(x) \rightarrow \int f(x) d\mu(x)\,.$$
  \end{itemize}
\end{prop}
The condition in {\it iv)} is the condition of tightness, or relative compactness, in $(M_1^p(\mathbb R), d_{W_p})$. In particular, it follows that,  for any $M \in \mathbb R$, the set
$$K_M := \{ \mu, \int |x|^p d\mu(x)  \leq M \}$$
is a compact set in $(M_1^q(\mathbb R), d_{W_q})$ for any $q <p$.

\section*{Acknowledgments}

We thank Gregory Schehr for the presentation of its model at the working group MEGA and Michel Ledoux for pointing out the semicontinuity in Remark \ref{Ledoux}.
\bibliographystyle{plain}

%\bibliography{schehr}

\end{document}